\documentclass[a4paper, 12pt]{amsart}

\usepackage{amsmath,amssymb,enumitem,verbatim,stmaryrd,xcolor,microtype,graphicx}
\usepackage[T1]{fontenc}
\usepackage[utf8]{inputenc}
\usepackage[english]{babel}
\usepackage[top=3.4cm,bottom=3.4cm,left=3.2cm,right=3.2cm]{geometry}
\usepackage[bookmarksdepth=3,linktoc=page,colorlinks,linkcolor={blue!80!black},citecolor={blue!80!black},urlcolor={blue!80!black}]{hyperref}

\usepackage{tikz}\usetikzlibrary{matrix,arrows,decorations.markings}
\usepackage{tikz-cd}
\pgfrealjobname{padicWanRiemannHyp}

\usepackage[mathcal]{euscript}                               
\usepackage{mathptmx}                                        
\usepackage{etoolbox}\makeatletter\patchcmd{\@startsection}{\@afterindenttrue}{\@afterindentfalse}{}{}\makeatother    
\patchcmd{\section}{\scshape}{\bfseries}{}{}\makeatletter\renewcommand{\@secnumfont}{\bfseries}\makeatother           
\usepackage[backgroundcolor=yellow,linecolor=yellow,textsize=footnotesize]{todonotes} \makeatletter \providecommand \@dotsep{5} \def\listtodoname{List of Todos} \def\listoftodos{\@starttoc{tdo}\listtodoname} \makeatother 

\theoremstyle{plain}
\newtheorem{thm}{Theorem}[section]
\newtheorem{cor}[thm]{Corollary}

\theoremstyle{definition}
\newtheorem{df}[thm]{Definition}
\newtheorem{rem}[thm]{Remark}
\newtheorem*{rem*}{Remark}

\newtheorem{conj}{Conjecture}

\usepackage{etoolbox}
\makeatletter
\patchcmd{\@startsection}{\@afterindenttrue}{\@afterindentfalse}{}{}             
\patchcmd{\part}{\bfseries}{\bfseries\LARGE}{}{}
\patchcmd{\section}{\scshape}{\bfseries}{}{}\renewcommand{\@secnumfont}{\bfseries} 
\patchcmd{\@settitle}{\uppercasenonmath\@title}{\large}{}{}
\patchcmd{\@setauthors}{\MakeUppercase}{}{}{}
\addto{\captionsenglish}{} 
\addto{\captionsenglish}{} 
\addto{\captionsenglish}{} 
\makeatother

\usepackage{fancyhdr}

\DeclareRobustCommand{\gobblefour}[4]{}    
\DeclareSymbolFont{sfoperators}{OT1}{bch}{m}{n} \DeclareSymbolFontAlphabet{\mathsf}{sfoperators} \makeatletter\def\operator@font{\mathgroup\symsfoperators}\makeatother 
\DeclareSymbolFont{cmletters}{OML}{cmm}{m}{it}              
\DeclareSymbolFont{cmsymbols}{OMS}{cmsy}{m}{n}
\DeclareSymbolFont{cmlargesymbols}{OMX}{cmex}{m}{n}
\DeclareMathSymbol{\myjmath}{\mathord}{cmletters}{"7C}     \let\jmath\myjmath 
\DeclareMathSymbol{\myamalg}{\mathbin}{cmsymbols}{"71}     
\DeclareMathSymbol{\mycoprod}{\mathop}{cmlargesymbols}{"60}
\DeclareMathSymbol{\myalpha}{\mathord}{cmletters}{"0B}     \let\alpha\myalpha 
\DeclareMathSymbol{\mybeta}{\mathord}{cmletters}{"0C}      \let\beta\mybeta
\DeclareMathSymbol{\mygamma}{\mathord}{cmletters}{"0D}     \let\gamma\mygamma
\DeclareMathSymbol{\mydelta}{\mathord}{cmletters}{"0E}     \let\delta\mydelta
\DeclareMathSymbol{\myepsilon}{\mathord}{cmletters}{"0F}   \let\epsilon\myepsilon
\DeclareMathSymbol{\myzeta}{\mathord}{cmletters}{"10}      \let\zeta\myzeta
\DeclareMathSymbol{\myeta}{\mathord}{cmletters}{"11}       \let\eta\myeta
\DeclareMathSymbol{\mytheta}{\mathord}{cmletters}{"12}     \let\theta\mytheta
\DeclareMathSymbol{\myiota}{\mathord}{cmletters}{"13}      \let\iota\myiota
\DeclareMathSymbol{\mykappa}{\mathord}{cmletters}{"14}     \let\kappa\mykappa
\DeclareMathSymbol{\mylambda}{\mathord}{cmletters}{"15}    \let\lambda\mylambda
\DeclareMathSymbol{\mymu}{\mathord}{cmletters}{"16}        \let\mu\mymu
\DeclareMathSymbol{\mynu}{\mathord}{cmletters}{"17}        \let\nu\mynu
\DeclareMathSymbol{\myxi}{\mathord}{cmletters}{"18}        \let\xi\myxi
\DeclareMathSymbol{\mypi}{\mathord}{cmletters}{"19}        \let\pi\mypi
\DeclareMathSymbol{\myrho}{\mathord}{cmletters}{"1A}       \let\rho\myrho
\DeclareMathSymbol{\mysigma}{\mathord}{cmletters}{"1B}     \let\sigma\mysigma
\DeclareMathSymbol{\mytau}{\mathord}{cmletters}{"1C}       \let\tau\mytau
\DeclareMathSymbol{\myupsilon}{\mathord}{cmletters}{"1D}   \let\upsilon\myupsilon
\DeclareMathSymbol{\myphi}{\mathord}{cmletters}{"1E}       \let\phi\myphi
\DeclareMathSymbol{\mychi}{\mathord}{cmletters}{"1F}       \let\chi\mychi
\DeclareMathSymbol{\mypsi}{\mathord}{cmletters}{"20}       \let\psi\mypsi
\DeclareMathSymbol{\myomega}{\mathord}{cmletters}{"21}     \let\omega\myomega
\DeclareMathSymbol{\myvarepsilon}{\mathord}{cmletters}{"22}\let\varepsilon\myvarepsilon
\DeclareMathSymbol{\myvartheta}{\mathord}{cmletters}{"23}  \let\vartheta\myvartheta
\DeclareMathSymbol{\myvarpi}{\mathord}{cmletters}{"24}     \let\varpi\myvarpi
\DeclareMathSymbol{\myvarrho}{\mathord}{cmletters}{"25}    \let\varrho\myvarrho
\DeclareMathSymbol{\myvarsigma}{\mathord}{cmletters}{"26}  \let\varsigma\myvarsigma
\DeclareMathSymbol{\myvarphi}{\mathord}{cmletters}{"27}    \let\varphi\myvarphi

\newcommand\C{{\mathbb C}}
\newcommand\F{{\mathbb F}}

\newcommand\Q{{\mathbb Q}}

\newcommand\Z{{\mathbb Z}}

\renewcommand\int{\textup{int}}

\renewcommand\geq{\geqslant}
\renewcommand\leq{\leqslant}

\newcommand{\DistTo}{\xrightarrow{
   \,\smash{\raisebox{-0.65ex}{\ensuremath{\scriptstyle\sim}}}\,}}

\title{$p$-adic Wan-Riemann Hypothesis for $\Z_p$-towers of curves}

\author{Roberto Alvarenga}
\address{\rm Instituto de Ci\^encias Matem\'aticas e Computac\~ao ICMC/USP, S\~ao Carlos, Brazil}
\email{{alvarenga@icmc.usp.br}}
\thanks{The author thanks the  Department of Mathematics of the University of California-Irvine (UCI) for hosting him during the preparation of this manuscript.}

\begin{document}

\begin{abstract}
Our goal in this paper is to investigate four conjectures, proposed by Daqing Wan,  about the stable behavior of a geometric $\Z_p$-tower of curves $X_{\infty}/X$. Let $h_n$ be the class number of the $n$-th layer in  $X_{\infty}/X$. It is known from Iwasawa theory that there are integers $\mu(X_{\infty}/X), \lambda(X_{\infty}/X)$ and $ \nu(X_{\infty}/X)$ such that the $p$-adic valuation $v_p(h_n)$ equals to $\mu(X_{\infty}/X) p^n + \lambda(X_{\infty}/X) n+ \nu(X_{\infty}/X)$ for $n$ sufficiently large.  Let $\Q_{p,n}$ be the splitting field (over $\Q_p$) of the zeta-function of $n$-th layer in $X_{\infty}/X$. The $p$-adic Wan-Riemann Hypothesis conjectures that  the extension degree  $[\Q_{p,n}:\Q_p]$ goes to infinity as $n$ goes to infinity. After motivating and introducing the conjectures, we prove the $p$-adic Wan-Riemann Hypothesis when $\lambda(X_{\infty}/X)$ is nonzero. 
\end{abstract}

\maketitle




\section{Introduction}
\label{introduction}


\subsection*{Intention and scope of this text}
Let $p$ be a prime number and $\Z_p$ be the $p$-adic integers. Our exposition is meant as a reader friendly introduction to the $p$-adic Wan-Riemann Hypothesis for geometric  $\Z_p$-towers of curves. Actually,  we shall  investigate four conjectures stated by Daqing Wan concerning the stable behavior of  zeta-functions attached  to a  geometric $\Z_p$-tower of curves, where these curves are defined over a finite field $\F_q$ and $q$ is a power of the prime $p$. The paper is organized as follows: in this first section we introduce the basic background of this subject  and motivate the study of $\Z_p$-towers;  in the Section 2, we state the four conjectures proposed by Daqing Wan (in particular the $p$-adic Wan-Riemann Hypothesis for geometric $\Z_p$-towers of curves) and re-write those conjectures in terms of some $L$-functions; in the last section, we prove the main theorem of this article using the $T$-adic L-function.

\subsection*{Motivating question} \label{motivatingquestion} Let $X$ be a smooth projective geometric irreducible curve defined over $\F_q$, and $N_n $ denotes the number of rational points of $X$ over a constant field extension $\F_{q^n}$ of $\F_q$. Consider
\[ Z(X,s) := \exp\left( \sum_{n \geq 0} \frac{N_n s^n}{n} \right)  = \prod_{x \in |X|} \frac{1}{1-s^{|x|}}\]
to be the zeta-function of $X$, where $|X|$ denotes the set of closed points of $X$ and $|x|$ the degree of the closed point $x$. The Riemann-Roch theorem implies that the zeta-function $Z(X,T)$ is a rational function in $s$ of the form  
\[Z(X,s) = \frac{P(X,s)}{(1-s)(1-qs)} \in \Q(s) \]
where $P(X,s) \in 1+s \Z[s]$ is a polynomial of degree $2g$, and $g$ is the genus of $X$ (see \cite[Theorem 6.1]{dino}).   
By a celebrated theorem of Weil \cite{Weil}, the analogous of Riemann Hypothesis for curves over finite fields (see \cite[Paragraphs 5.8 and 5.9, pag.\ 283]{dino}),  
\[
P(X,s) := \prod_{i=1}^{2g} (1-\alpha_i s) \in \C(s)
\]
 where $\alpha_i$ are algebraic integers with complex  norm $|\alpha_i| = \sqrt{q},\ i=1, \ldots,2g.$

The following natural question arises in this context. Let $\Q_p$ be the field of $p$-adic numbers and $\C_p$ be the completion of $\overline{\Q}_p$, the algebraic closure of $\Q_p.$ The fields $\C$ and $\C_p$ are isomorphic as abstract fields (see \ \cite[III.3.5]{robert}). Considering 
$$P(X_,s) = \prod_{i=1}^{2g_n} (1-\alpha_i s) \in \C_p(s)$$
 as a decomposition of $P(X,s)$ over $\C_p$,  we may ask about $|\alpha_i|_p := p^{-v_p(\alpha_i)}$ the $p$-adic norm  of $\alpha_i,\ i=1, \ldots, 2g.$ An answer for this question could be stated as a $p$-adic version of the classical Riemann Hypothesis. 
The $p$-adic valuation $v_p(\alpha_i)$, also called \textbf{$p$-slopes} of $X$, remains quite mysterious in general  and does not seem to have a stable behavior as the complex norm $|\alpha_i|.$ Hence, the $p$-slopes of $X$ have been studied through the Newton polygon $NP(X)$ of $X$, which is by definition the convex hull of the points $(i, v_p(\alpha_i)), i=1, \ldots, 2g$.  By observing the Weil theorem we can order the reciprocal roots of $P(X,s)$ such that 
\[
0 \leq v_p(\alpha_1) \leq \cdots \leq v_p(\alpha_{2g}) \leq r
\]
where $r$ is such that $q=p^r.$  

Another approach to studying the behavioral problem of $p$-slopes of $X$, is to
understand these rational numbers when $X$ varies, for instance, in an algebraic
family or when $X$ varies in a $Z_p$-tower (which will be our specific interest here, in the spirit of
Iwasawa theory).

\subsection*{The $\Z_p$-towers }
 In the late 1950s, Kenkichi Iwasawa proved several results and formulated some important conjectures concerning the behavior of the ideal class groups in a certain $\Z_p$-tower of number fields, see \cite{iwasawa1} and \cite{iwasawa2}.  After Iwasawa, some authors started developing his ideas not just for number fields, but also for global function fields.  This  is what is now known as  Iwasawa theory. 
By definition, a $\Z_p$-tower (or $\Z_p$ extension) of a field $F$ is an infinite Galois extension $F_{\infty}/F$  with Galois group isomorphic to the additive group of $p$-adic integers $\Z_p.$ Since the closed subgroups of $\Z_p$ are of the form $0$ and $p^n\Z_p$, it is possible to regard a $\Z_p$-tower $F_{\infty}/F$  as a sequence of fields
\[ F=F_0 \subseteq F_1 \subseteq \cdots \subseteq F_{\infty}  \]
with Galois groups $\mathrm{Gal}(F_n,F) \simeq \Z/p^n \Z$ and $F_{\infty} = \bigcup_{n=1}^{\infty} F_n$. We say $F_n$ is the $n$-th layer of the $\Z_p$-tower.

Given a global function field $F$ over the finite field $\F_q$, the constant field extension $ \F_{q^{p^{\infty}}} F/F$ is an immediate $\Z_p$-tower.  By a \textbf{geometric extension} of $F$, we mean an extension which does not contain any constant field extension. A \textbf{geometric $\Z_p$-tower}  $F_{\infty}/F$ is the one  such that $F_n/F$ is a geometric extension for every $n$-th layer $F_n.$ 

Let $F_{\infty}/F$ be a constant $\Z_p$-tower, i.e.\ $F_{\infty} =\F_{q^{p^{\infty}}} F$ and $F$  a global function field over $\F_q$, then there are no ramified primes in $F_{\infty}$ (\cite[Prop. 8.5]{rosen}).
For a given $\Z_p$-tower of number fields, there are only finitely many ramified primes.  In contrast,  there may be infinitely many ramified primes in a geometric $\Z_p$-tower (see \cite{gold}). Throughout of this paper, we shall just work with geometric $\Z_p$-towers with finitely many ramified primes. The condition is natural and necessary in order for many
definitions and theorems to make sense.

\subsection*{Iwasawa class number formula} Assume $F_{\infty}/F$ be a $\Z_p$-tower over a global field $F$. Let $h_n$ denote the class number of $n$-th layer $F_n$ in  $F_{\infty}/F.$   It is well known from Iwasawa theory that there are integers $\mu(F_{\infty}/F), \lambda(  F_{\infty}/F)$ and $\nu(  F_{\infty}/F)$ such that 
\begin{equation} \label{iwsawaformula} 
v_p(h_n) =\mu(F_{\infty}/F)p^n + \lambda(  F_{\infty}/F) n+ \nu(  F_{\infty}/F) 
\end{equation}
for all sufficiently large $n$. This formula was proved by Iwasawa  first when $ F_{\infty}/F$ is a constant $\Z_p$-tower of  global function fields, and thus   for $  F_{\infty}/F$  a $\Z_p$-tower  of number fields, see \cite[Thm. 4]{iwasawa1}. Later, was noticed by Mazur-Wiles \cite{wiles}  and Gold–Kisilevsky \cite{gold} that Iwasawa's idea also works for geometric $\Z_p$-towers  of global function fields.  We could also wonder about the behavior of $h_n$ in a $\Z_{p}^d$-tower, for any $d \geq 1$, see \cite{wan-class} for a recent work in this direction.

\subsection*{Acknowledgments} The author wishes to express his gratitude to Daqing Wan, who teaches him almost everything that he knows about the subject of this article. He also would like to thank Joe Kramer-Miller for many fruitful conversations. The author was supported by FAPESP grant numbers [2017/21259-3] and [2018/18805-9].


\section{The Wan conjectures}

By a $\Z_p$\textbf{-tower of curves} we mean the $\Z_p$-tower of its function fields. Consider a $\Z_p$-tower 
\[
X_{\infty}: \cdots \rightarrow X_n \rightarrow \cdots \rightarrow X_1 \rightarrow X_0 = X
\]
of smooth projective geometric irreducible curves defined over $\F_q$. As before, the above $\Z_p$-tower gives continuous isomorphism
\[
\rho: G_{\infty} := \mathrm{Gal}(X_{\infty}, X) \DistTo\Z_p,
 \] 
and for each $n$, 
\[ 
G_n := \mathrm{Gal}(X_{n}, X) \simeq \Z \big/p^n\Z.
\]
Let $S$ be the ramification locus of the tower, which is a subset of $|X|$ the closed points of $X.$ The tower is unramified on its complement $U := |X| - S.$ 
By class field theory, the ramification locus $S$ is non-empty and for each non-empty finite $S$, there are uncountable many such $\Z_p$-towers over $\F_q.$ 
  We shall assume throughout the paper that $S$ is finite. For each $n \geq 0$, let
\[
Z(X_n,s) := \prod_{x \in |X_n|} \frac{1}{1-s^{|x|}} \in 1 + s \Z[[s]]
\]
be the zeta-function of $X_n$, where $|x|$ denotes the degree of $x$. As above, we may write

\[ 
Z(X_n,s) = \frac{P(X_n,s)}{(1-s)(1-qs)} \in \Q(T)
\]
where $P(X_n,s) \in 1+s \Z[s]$ is a polynomial of degree $2g_n$, and $g_n$ is the genus of $X_n$.
Considering the decomposition  
\[
P(X_n,s) := \prod_{i=1}^{2g_n} (1-\alpha_i(n)s) \in \C_p[s],
\]
our motivating question (Section \ref{motivatingquestion}) arises in this context as follows: how $v_p(\alpha_i(n))$ varies, for $i=1, \ldots, 2g_n$, as $n$ goes to infinity? 
We shall investigate this problem in the aspect of the splitting field of $P(X_n,s)$ over $\Q_p$, as $n$ goes to infinity.

We are now in position to state the Daqing Wan conjectures (\cite[Conj. 3.3]{wan-zptowers}) concerning the degree extension of the splitting field of $P(X_n,s)$ over $\Q_p$ (and $\Q$).   

\begin{conj} \label{wanconj1}Let $\Q_n$  denote the splitting field of $P(X_n,s)$ over $\Q$, namely  $\Q_n := \Q(\alpha_i(n), i=1,\ldots, 2g_n ).$ 
Then the extension degree $[\Q_n:\Q]$ goes to infinity as $n$ goes to infinity. 

\end{conj}

\begin{conj}[\textbf{$p$-adic Wan-Riemann Hypothesis}] \label{wanconj2}Let $\Q_{p,n}$ denote  the splitting field of $P(X_n,s)$ over $\Q_p$, namely $\Q_{p,n} := \Q_p(\alpha_i(n), i=1,\ldots, 2g_n ).$ Then the extension degree $[\Q_{p,n}:\Q_p]$ goes to infinity as $n$ goes to infinity. 
\end{conj}

\begin{conj}\label{wanconj3}  Let $\Q_{p,n}$ denote  the splitting field of $P(X_n,s)$ over $\Q_p$, namely $\Q_{p,n} := \Q_p(\alpha_i(n), i=1,\ldots, 2g_n ).$ Then  the ramification degree $[\Q_{p,n}:\Q_p]^{\mathrm{ram}}$ goes to infinity as $n$ goes to infinity. 

\end{conj}

\begin{conj} \label{wanconj4}  Let $\Q_{p,n}$ denote  the splitting field of $P(X_n,s)$ over $\Q_p$, namely $\Q_{p,n} := \Q_p(\alpha_i(n), i=1,\ldots, 2g_n ).$ Then  there is a positive constant $c(X_{\infty}/X)$ depending on the $\Z_p$-tower $X_{\infty}/X$ such that for all sufficiently large $n$, $[\Q_{p,n}:\Q_p]^{\mathrm{ram}} \geq c(X_{\infty}/X) p^n.$
\end{conj}

\begin{rem} \label{remark-conj} Observe that each conjecture is stronger than the previous one. That Conjecture \ref{wanconj4} implies Conjecture \ref{wanconj3} and  Conjecture \ref{wanconj3} implies Conjecture \ref{wanconj2} is immediate. In order to show how the $p$-adic Wan-Riemann Hypothesis implies Conjecture \ref{wanconj1}, see \cite[Exerc. 5, pag. 74]{koblitz}. 

\end{rem}

\subsection*{Impact on motivating question} The Conjecture \ref{wanconj3} says that, as  $n$ grow, there is a reciprocal root $\alpha_i(n)$ of $P(X_n,s)$ such that $v_p(\alpha_i(n))$ goes to zero. In terms of the Newton polygon of $X_n$, it says that the first slope of $NP(X_n)$ goes to zero as $n$ goes to infinity.   

  \medskip

In \cite[Thm. 4]{gold}, Gold and Kisilevsky prove as a corollary of Theorem 4, that Conjecture \ref{wanconj1} holds when $\mu(F_{\infty}/F)=0$ and $\lambda(F_{\infty}/F) \neq 0$ in the Iwasawa Formula (\ref{iwsawaformula}).  Our main contribution to this problem, Theorem \ref{maintheorem} in the next section, is to prove Conjecture \ref{wanconj4} (and thus Conjectures \ref{wanconj3}, \ref{wanconj2} and \ref{wanconj1}) when $\lambda(F_{\infty}/F) \neq 0.$

Before reframe the Wan conjectures in terms of certain $L$-functions, let us explain our motivation to call  Conjecture \ref{wanconj2} of $p$-adic Wan-Riemann Hypothesis. 
In \cite{wan-douglas}  Wan and Haessig continue the investigation initiated by Wan in \cite{wan-cycles} about the zeta-function $Z_r(Y,s)$ of $r$-cycles in a projective variety $Y$ defined over a finite field. They conjecture in that article that the splitting field of $Z_r(Y,s)$ over $\Q_p$ is a finite extension of $\Q_p.$  This conjecture is the \cite[Conj. III]{wan-douglas} and is called $p$-adic Riemann Hypothesis. They were motivated by Goss \cite{goss} (where the author re-examines Wan's work \cite{wan-rh}) and his corresponding formulation for the characteristic $p$ $L$-function.  Among others things, Wan and Haessig  prove this conjecture when  $Y$ is a normal connected scheme of dimension $\dim Y \geq 2$, $r = \dim Y-1$, and the Chow group of $r$-cycles is rank one. 


\subsection*{$L$-functions and Wan's conjectures}  We finish this section rewriting the above conjectures in terms of classical $L$-functions. First, we observe the following consecutive divisibility
$ P(X_n,s) \big|  P(X_{n+1},s),$ for $\ n \geq 0.$
Hence, in order to investigate how changes the extension degrees in the Conjectures \ref{wanconj1}, \ref{wanconj2},  \ref{wanconj3} and \ref{wanconj4}, as $n$ grow, it is enough to study for $n \geq 1$, the primitive part of $Z(X_n,s)$ defined by
\[
Q(X_n,s) := \frac{ P(X_n,s)}{ P(X_{n-1},s)} = \frac{ Z(X_n,s)}{ Z(X_{n-1},s)}. 
\] 

Since the Galois group $\mathrm{Gal}(X_m,X_n)$ is a cyclic $p$-group for all $m \geq n$, if $X_n$ is ramified over $X$ at a closed point $x \in S,$ then $X_m$ is totally ramified over $X_n$ at $x$. Therefore, without loss of generality, by going to a larger $n$ if necessary, we can assume that $X_1$ is already ramified at every point of $S$. From now on, we assume that $X_1$ is indeed (totally) ramified at every point of $S.$ 

For a primitive character  $\chi_n : G_n \rightarrow \C_{p}^{\ast}$
of order $p^n > 1,$ we consider the $L$-function of $\chi_n$ over $X$ 
\[
L(\chi_n, s) = \prod_{x \in |U|} \frac{1}{1- \chi_n(\mathrm{Frob}_x) s^{|x|}} \  \  \in 1 + s \Z[\chi_n(1)][[s]],
\]
where $|U|$ denotes the set of closed points of $U$ and $\mathrm{Frob}_x$ denotes the arithmetic Frobenius element of $G_n$ at $x.$ 
Besides that, since $\chi_n$ is primitive of order $p^n$, $\zeta_{p^n} := \chi_n(1)$ is a primitive $p^n$-th root of unity. 

By a theorem due to Weil (\cite[Thm. 9.16B]{rosen}), the $L$-function $L(\chi_n, s)$ is a polynomial in $s$ of degree $\ell(n)$ and its reciprocal roots have complex norm equals to $\sqrt{q}.$ Moreover, 
we have the decomposition
\[
Q(X_n,s) = \prod_{\chi_n: G_n \rightarrow \C_{p}^{\ast}} L(\chi_n,s),
\] 
where the product takes all $\chi_n$ primitive characters of order $p^n,$ see \cite[Paragraph 4, pag.130]{rosen} for more details.  

Let $\sigma \in \mathrm{Gal}(\Q(\zeta_{p^n}),\Q) =\mathrm{Gal}(\Q_p(\zeta_{p^n}),\Q_p),$ one checks that 
\[
L(\chi_n,s)^{\sigma} = L(\chi_{n}^{\sigma},s).
\]
It follows that the degree and the slopes of $L(\chi_n,s)$ depend only on $n,$ not on the choice of the primitive character $\chi_n$ of $G_n.$ Therefore, in order to study the extension degree conjectures for the splitting field of $P(X_n,s)$, we can just choose and fix one primitive character $\chi_n$ of order $p^n$ for each $n \geq1.$  Hence, the  Conjectures \ref{wanconj2}, \ref{wanconj3} and \ref{wanconj4} are reduced to: 

\begin{conj} \label{conj-wan2.1} Let $L_{p,n}$ denote the splitting field of $L(\chi_n, s)$ over $\Q_p.$ 
The extension degree $[L_{p,n}:\Q_p]$ goes to infinity as $n$ goes to infinity. 

\end{conj}

\begin{conj} \label{conj-wan3.1} Let $L_{p,n}$ denote the splitting field of $L(\chi_n, s)$ over $\Q_p.$ 
The ramification degree $[L_{p,n}:\Q_p]^{\mathrm{ram}}$ goes to infinity as $n$ goes to infinity. 

\end{conj}

\begin{conj} \label{conj-wan4.1} Let $L_{p,n}$ denote the splitting field of $L(\chi_n, s)$ over $\Q_p.$ 
Then  there is a positive constant $c(X_{\infty}/X)$ depending on the $\Z_p$-tower $X_{\infty}/X$ such that for all sufficiently large $n$, $[L_{p,n}:\Q_p]^{\mathrm{ram}} \geq c(X_{\infty}/X) p^n.$
\end{conj}


\section{$T$-adic L-function and the main theorem}

In order to prove the $p$-adic Wan-Riemann Hypothesis for geometric $\Z_p$-tower of curves, we need to introduce the $T$-adic $L$-function (first defined  in \cite{liu-wan}) associated to a $\Z_p$-tower of curves. 

In the previous section we have rewritten Wan's conjectures by considering the $L$-function associated to a finite character $\chi_n :G_n \rightarrow \C_{p}^{\ast}$. Now, we also consider all continuous $p$-adic characters $\chi: G_{\infty} \rightarrow \C_{p}^{\ast},$ not necessarily of finite order. 

The isomorphism 
\[ 
\rho: G_{\infty}   \DistTo\Z_p,
\]
given by the geometric $\Z_p$-tower $X_{\infty}/X$ of smooth projective geometrically irreducible curves over $\F_q$, induces
\begin{align*}
f_{\rho}: |U| &\longrightarrow  \Z_p \\
x &\longmapsto \rho(\mathrm{Frob}_x) 
\end{align*}
the $p$-adic valued Frobenius function. This function determines the $\Z_p$-tower $X_{\infty}/X$ by class field theory, i.e.\ any condition we would impose on the tower is a condition  on this Frobenius function $f_{\rho}.$ 
Consider the universal continuous $T$-adic character 
\begin{align*}
\Z_p & \longrightarrow  \Z_p[[T]]^{\ast} \\
1 & \longmapsto  1+T . 
\end{align*}
Composing this universal $T$-adic character of $\Z_p$ with the isomorphism $\rho$, we have the universal $T$-adic character of $G_{\infty}$
\[
\rho_T : G_{\infty} \rightarrow \Z_p \rightarrow \Z_p[[T]]^{\ast}.
\]
Let $D_p(1)$ denote the open unit disk in $\C_p.$ For any element $t \in D_p(1)$, we have the following evaluating map 
\begin{align*}
\Z_p[[T]]^{\ast} & \longrightarrow \C_{p}^{\ast} \\
T & \longmapsto  t.
\end{align*}
Composing all these maps, we obtain, for a fixed $t \in D_{p}(1)$, a continuous character 
\begin{eqnarray}
\rho_t : G_{\infty} \longrightarrow \C_{p}^{\ast}.
\end{eqnarray}
The $L$-function of $\rho_t$ is defined by
\[
L(\rho_t,s) := \prod_{x \in |U|} \frac{1}{1-\rho_t(\mathrm{Frob}_s) s^{|x|}} \  \  \in 1+ s\C_p[[s]].
\]

\begin{rem} The open unit disk $D_p(1)$ parameterizes all continuous $\C_p$-valued characters $\chi$ of $G_{\infty}$ via the relation
$ t=\chi(1) -1.$
In the case that $\chi = \chi_n$ is a $p$-adic character of $G_{\infty}$ of order $p^n,$ then $\chi(1)$ is a primitive $p^n$-th root of unity. Hence, 
$t = t_n := \chi_n(1) -1,$ and  $ | t_n |_p = p^{- \frac{1}{p^{n-1}(p-1)}}.$ Therefore, 
\[
L(\chi_n,s) = L(\rho_{t_n},s).
\]
\end{rem}

Elements of the form $t_n := \chi_n(1)-1$, for $n \geq 0$, are called the \textbf{classical points} in $D_p(1).$ As $n$ goes to infinity, $t_n$ approaches to the boundary of the disk $D_p(1).$ Thus, in order to understand the behavior of $L(\chi_n,s)$ as $n$ grows, in particular Conjectures \ref{conj-wan2.1}, \ref{conj-wan3.1} and \ref{conj-wan4.1}, it is enough to understand the $L$-function $L(\rho_t,s)$ for all $t$ near the boundary of $D_p(1).$ Namely, we should consider the following Liu-Wan universal $L$-function. 

\begin{df} The \textbf{$T$-adic $L$-function} $L_{\rho}(T,s)$ of the tower $X_{\infty}/X$ is by definition the $L$-function of the $T$-adic character $\rho_T$
\[
L_{\rho}(T,s) := L(\rho_T,s) = \prod_{x \in |U|} \frac{1}{1-(1+T)^{\rho(\mathrm{Frob}_x)} s^{|x|}} \  \ \in 1+s\Z_p[[T]][[s]].
\]
This is a $p$-adic power series in the two variables $T$ and $s.$ 
\end{df}

For $t \in D_p(1)$, we have 
\[
L_{\rho}(T,s)|_{T=t} = L_{\rho}(t,s) = L(\rho_t,s).
\]
 If $t=t_n$ is a classical point, $L_{\rho}(T,s)|_{T=t_n} =L(\chi_n,s)$, as noted above.   
We can now formulate and prove our main theorem. 

\begin{thm} \label{maintheorem} Let $X_{\infty}/X$ be a geometric $\Z_p$-tower of smooth projective geometrically irreducible curves over a finite field $\F_q.$ Suppose $X_{\infty}/X$ ramified at a finite subset $S$ of $|X|$. If $\lambda(X_{\infty}/X)$ is nonzero, then  Conjecture \ref{wanconj4} is true. 

\begin{proof} Replacing $X$ by $X_m$, for some $m \geq 1$ if necessary, we can assume $X_{\infty}/X$ is totally ramified at the finite subset $S$ of $|X|$, and unramified at $U = |X| - S.$ 

We first observe   (following \cite[Thm. 3.5]{wan-zptowers}) that we can write the specialization at $s=1$, of the $T$-adic $L$-function $L_{\rho}(T,s)$ of $X_{\infty}/X$ as following
\[
L_{\rho}(T,1) = p^{\mu} (T^{\lambda} + p a_1 T^{\lambda -1} + \cdots + p a_{\lambda})u(T),
\]
where $a_i \in \Z_p,$ for $ i=1, \ldots, \lambda$, $u(T)$ is a unit in $\Z_p[[T]],$ and $\mu = \mu(X_{\infty}/X), \lambda = \lambda(X_{\infty}/X)$ as in the Iwasawa class number formula (\ref{iwsawaformula}). Write
\[
L_{\rho}(T,s) = \sum_{k=0}^{\infty} L_k(T)s^k,
\]
for $L_k(T) \in \Z_p[[T]].$ Since $L_{\rho}(t_n,s)$ is a polynomial of degree $\ell(n)$, $L_k(t_n)=0$ for all $k > \ell(n).$ The $p$-adic Weierstrass preparation theorem (see \cite[Thm 14, pag. 105]{koblitz}) implies that 
\[
L_{k}(T) = \frac{(1+T)^{p^n} - 1}{T} u_k(T),
\]
where $u_k(T) \in \Z_p[[T]]$, for all $k> \ell(n).$ It follows that the series $L_{\rho}(T,s)$ is $(p,T)$-adically convergent for $s \in \Z_p[[T]].$ If $s=1$, note that $L_{\rho}(T,1) \neq 0$ as its specialization at classical points $t_n$ is nonzero. Therefore we can write
\[
 L_{\rho}(T,1) = p^{\mu} (T^{\lambda} + p a_1 T^{\lambda -1} + \cdots + p a_{\lambda})u(T),
 \]
with $a_i \in \Z_p$, for $ i=1, \ldots, \lambda$, and $u(T)$ is a unit in $\Z_p[[T]].$ We are left to show that $\mu$ and $\lambda$ are the Iwasawa constants (\ref{iwsawaformula}) for the $\Z_p$-tower $X_{\infty}/X.$  Let $h_n$ (resp. $h$) be the class number of the $n$-th layer $X_n$ (resp. $X$). Observe that
\[
v_p(h_n) - v_p(h) = \sum_{k=1}^{n} p^{k-1}(p-1) v_p(L_{\rho}(t_k,1)).
\]
From the above identity, 
$
L_{\rho}(t_k,1) = p^{\mu} (t_{k}^{\lambda} + p a_1 t_{k}^{\lambda -1} + \cdots + p a_{\lambda})u(t_k).
$
Since $u(T)$ is a unit in $\Z_p[[T]],$ $v_p(u(t_k)) =0.$ Thus, for $p^{k-1}(p-1) > \lambda$, we have
\[
v_p(t_{k}^{\lambda} + p a_1 t_{k}^{\lambda -1} + \cdots + p a_{\lambda}) = \frac{\lambda}{p^{k-1}(p-1)},
\]
since $v_p(t_k) = 1/p^{k-1}(p-1).$  For $n$ sufficiently large we conclude that 
\[
v_p(h_n) = \mu p^n + \lambda n + \nu
\]
for some constant $\nu$. Therefore, $\mu = \mu(X_{\infty}/X)$ and $\lambda = \lambda(X_{\infty}/X).$

Next we consider $\chi_n : \Z/p^n\Z \rightarrow \C_{p}^{\ast}$ a primitive character of order $p^n > 1$, for $n \geq 1.$   As before, we can write the $L$-function of $\chi_n$ over $X$ as follows
\[
L(\chi_n,s) = \prod_{x \in |U|} \frac{1}{1- \chi_n(\mathrm{Frob}_x) s^{|x|}} = \prod_{i=1}^{\ell(n)} (1- \beta_i(n)s) \in \C_p[s].
\]
Since
\[
Q(X_n,s) = \prod_{\chi_n : G_n \rightarrow \C_{p}^{\ast}} L(\chi_n,s) ,
\]
where the product is taken over all primitive character of order $p^n,$ we see that $L(\chi_n,s)$ divides $P(X_n,s).$  Thus,
\[
\Q_{p,n} \supseteq \big\{ \alpha_i(n) \big| 1 \leq i \leq 2g_n \big\} \supseteq \big\{ \beta_i(n) \big| 1 \leq i \leq \ell(n) \big\},
\] 
where as in the Conjecture \ref{wanconj4},  $\alpha_i(n)$ are the reciprocal roots of $P(X_n,s)$, for $i=1, \ldots, 2g_n$, and $\Q_{p,n}$ is the splitting field of $P(X_n,s)$ over $\Q_p$. Hence
\[
L(\chi_n,1) = \prod_{i=1}^{\ell(n)} (1- \beta_i(n)) \in \Q_{p,n}.
\]
Moreover, we have already observed that $L(\chi_n,1) = L_{\rho}(t_n,1)$, where $t_n = \chi_n(1)-1.$ It follows from the above description of $ L_{\rho}(T,1)$ that
\begin{eqnarray*}
v_p(L(\chi_n,1)) &=& \mu + v_p(t_{n}^{\lambda} + p a_1 t_{n}^{\lambda -1} + \cdots + p a_{\lambda}) \\
&=& \mu + \frac{\lambda}{p^{n-1}(p-1)}
\end{eqnarray*}
for $p^{n-1}(p-1) > \lambda.$ Since by hypothesis $\lambda \neq 0$, 
\[
v_p(L(\chi_n,1)) = \mu + \frac{\lambda}{p^{n-1}(p-1)},
\]
for all $n$ such that $p^{n-1}(p-1) > \lambda.$ Therefore, 
\[
[\Q_{p,n}:\Q_p]^{\mathrm{ram}} \geq \frac{p^{n-1}(p-1)}{\gcd(p^{n-1}(p-1),\lambda)} \geq \frac{p^{n-1}(p-1)}{\lambda} = c p^n,
\]
where $c = c(X_{\infty}/X) =(1-p^{-1})\lambda^{-1} > 0.$ That is,  Conjecture \ref{wanconj4} is true.   
   \end{proof}
\end{thm}

From Remark \ref{remark-conj}, the above theorem has the $p$-adic Wan-Riemann Hypothesis as corollary  when $\lambda(X_{\infty}/X) \neq 0$.  

\begin{cor} Let $X_{\infty}/X$ be a geometric $\Z_p$-tower of smooth projective geometrically irreducible curves over a finite field $\F_q.$ Suppose $X_{\infty}/X$ ramified at a finite subset $S$ of $|X|$. If $\lambda(X_{\infty}/X)$ is nonzero. Then Conjecture \ref{wanconj3}  is true.  
$\hfill \square$

\end{cor} 

\begin{cor} Let $X_{\infty}/X$ be a geometric $\Z_p$-tower of smooth projective geometrically irreducible curves over a finite field $\F_q.$ Suppose $X_{\infty}/X$ ramified at a finite subset $S$ of $|X|$. If $\lambda(X_{\infty}/X)$ is nonzero. Then the $p$-adic Wan-Riemann Hypothesis (Conjecture \ref{wanconj2}) is true.  
$\hfill \square$

\end{cor}

\begin{cor} Let $X_{\infty}/X$ be a geometric $\Z_p$-tower of smooth projective geometrically irreducible curves over a finite field $\F_q.$ Suppose $X_{\infty}/X$ ramified at a finite subset $S$ of $|X|$. If $\lambda(X_{\infty}/X)$ is nonzero. Then  Conjecture \ref{wanconj1} is true.  
$\hfill \square$
\end{cor}

\medskip


\end{document}